\documentclass[11pt, fleqn]{amsart}
\usepackage{amssymb,amsmath,latexsym,amsbsy,color,enumerate}
\numberwithin{equation}{section}
\def\NN{\mbox{$I\hspace{-.06in}N$}}

\def\RR{\mbox{$I\hspace{-.06in}R$}}
\def\CC{\mbox{$C\hspace{-.11in}\protect\raisebox{.5ex}{\tiny$/$}
\hspace{.06in}$}}

\def\h{\hspace*{.24in}}

\newtheorem{theorem}{Theorem}
\newtheorem{definition}[theorem]{Definition}
\newtheorem{lemma}{Lemma}
\newtheorem{remark}{Remark}
\newtheorem{proposition}{Proposition}

\begin{document}
 \title{Stability and Identity of analytic functions in Hardy spaces}
    \author{Dang Duc Trong}
    \address{Dang Duc Trong, Department of Mathematics, Hochiminh City National University, 227 Nguyen Van Cu, Q5, HoChiMinh City, Vietnam}
    \email{ddtrong@mathdep.hcmuns.edu.vn}
\author{Truong Trung Tuyen}
    \address{Truong Trung Tuyen, Department of Mathematics, Hochiminh City National University, 227 Nguyen Van Cu, Q5, HoChiMinh City, Vietnam}
 \email{truongtrungtuyen2005@yahoo.com}
    \date{\today}
\maketitle
\section{Introduction}
Let $E$ be a subset of the unit disc $U$ of the complex plane $\CC$. Throughout this paper, we assume that $E$ consists of infinitively many points. We wish to reconstruct a function $f$ in the Hardy space $H^p(U)$, $1$ $\leq$ $p$ $\leq$ $\infty$, when its value at any point in $E$ is given. Recall that $H^p(U)$ is the space of all holomorphic functions $g$ on $U$ for which $\|g\|_p$ $<$ $\infty$, where
\begin{align*}
\|g\|_p &= \lim_{r\uparrow 1} \left\{\frac{1}{2\pi}\int_0^{2\pi}|g(re^{i\theta})|^p d\theta\right\}^{1/p}\h(1 \leq p < \infty),\\
\|g\|_{\infty} &= \lim_{r\uparrow 1} \sup_{\theta} |g(re^{i\theta})|.
\end{align*}
For convenience, from now on, we will denote $H^p(U)$ by $H^p$.

As it is well known, this is an ill-posed problem. Clearly, the interpolation $f$ is unique only if
\begin{description}
\item[(B)] $E$ contains a non-Blaschke sequence $(z_j)$, that is, a sequence satisfying the condition
\[
\displaystyle\sum_{j=1}^\infty(1-|z_j|) = \infty.
\]
\end{description}

Put
\begin{equation}
C_p(\varepsilon, R) = \sup \{\sup_{|z| \leq R}|g(z)|: g\in H^p,\, \|g\|_p\leq 1,\,|g(\zeta)| \leq \varepsilon\;\forall \zeta\in E\},
\label{CeRdefinition}
\end{equation}
for positive $\varepsilon$ and $R$ in $(0, 1)$. It can be seen that $C_p(\varepsilon, R)$ is bounded from above by $(1-R^2)^{-1/p}$. We consider the stability of the problem in the sense that $C_p(\varepsilon, R)$ converges to $0$ as $\varepsilon$ decreases to $0$. We also wish to establish the rate of convergence. In general, given only that $E$ satisfies (B), it is impossible to say anything about the rate of convergence. However, the convergence itself can be attained.

\begin{theorem}
If $E$ satisfies (B), then
\[
\lim_{\varepsilon \rightarrow 0}C_p(\varepsilon, R) = 0.
\]
\label{ConvergenceItselfTheorem}
\end{theorem}

\begin{remark}
When $E$ does not satisfy (B), the conclusion of Theorem \ref{ConvergenceItselfTheorem} is obviously false. As an example, consider the case $E$ is the sequence
$\{z_k\}$ where $\sum (1 - |z_k|)$ $<$ $\infty$. Let $g$ be the Blaschke product whose zeros are the points in $E$,
$g(z)$ $=$ $\displaystyle\prod \frac{z_k - z}{1 - \overline{z_k}z}\frac{|z_k|}{z_k}$. Then $g$ is in $H^{\infty}$, $\|g\|_{\infty}$ $=$ $1$, $g(\zeta)$ $=$ $0$
for all $\zeta$ in $E$, yet $g$ is not identically zero in $U$.
\end{remark}

To estimate the rate of convergence, we need some other materials, especially, the properties of the set $E$. For example, if $p$ $=$ $\infty$ and $E$ is a curve
of a simply connected domain $D$, one might use the harmonic measure $\omega$ of $E$ with respect to $D$ to get $|f(z)|$ $\leq$ $\varepsilon^{\omega(z)}$
inside $D$. Lately, in \cite{LRS}, Lavrent'ev, Romanov and Shishat-skii used a certain characteristic of the projection of $E$ onto the real axis. They showed that
 if $E$ $\subset$ $D$ $=$ $\{z: |z| \leq 1/4\}$ then $C_p(\varepsilon, R)$ $\leq$ $\max\{\varepsilon^{4/25},(6/7)^{n(\varepsilon)}\}$ for all $R$ $\in$ $(0,1/4)$,
  in which $n(\varepsilon)$ $\rightarrow$ $\infty$ as $\varepsilon$ $\rightarrow$ $0$. This approach is quite interesting that $E$ could be a sequence but must
  lye strictly inside $U$. Further, to our knowledge, all such estimates just give us the upper bounds. In this paper, we introduce some new quantities and then
  establish a quasi-polynomial two-side bound for $C_p(\varepsilon, R)$. We stress that the lower bound is also included in the estimate.
  Another advantage of this method is that the requirement that $\overline E$ lies inside $U$can be relaxed. In case $\overline E$ $\subset$ $U$, we get the
  following result
\begin{theorem} If $\overline E$ $\subset$ $U$ then there exists $\epsilon _0>0$ such that for $0<\epsilon <\epsilon _0$ there is correspondingly a finite Blaschke product $B_{\epsilon}(z)$ whose zeros are in $\overline E$ satisfying

\begin{eqnarray*}
\max _{|z|\leq R}|B_\epsilon (z)|\leq C_p(\epsilon ,R)\leq C\max _{|z|\leq R}|B_\epsilon (z)|^{1/2},
\end{eqnarray*}
where $C$ is a positive constant that depends only on $R$ and $p$. Moreover we have
\begin{eqnarray*}
\sup _{z\in E}|B_{\epsilon}(z)|\leq \epsilon.
\end{eqnarray*}
\label{LRSTheorem}
\end{theorem}

Closely related to the problem of stability of reconstructing analytic functions is the topic of identity of analytic functions. That is, under what conditions a bounded analytic function equals to zero identitically. It is a classical result that if $f$ is a function of class Nevalina and takes value zero on a non-Blaschke sequence, then $f$ is the zero function. It is an interseting question to ask is there a lower bound for the fastness of tending to zero of analytic functions $|f(z)|$ when $|z|\rightarrow 1$ along $E$ assuming that $E$ is a set having $\overline{E}\cap \partial U\not= \emptyset$? Again, this question is trivial if $E$ contains no non-Blaschke sequence. For the case $E$ is a non-Blaschke sequence (not neccesary have limit points in $\partial U$) containing no limit point of itself, this question was thoroughly explored in \cite{dan} and \cite{hay} (see also \cite{lk}).

In particular, Danikas\cite{dan} and Hayman\cite{hay} proved that if $(z_j)$ is a non-Blaschke sequence having no limit point of itself, there is a sequence of positive real numbers $(\eta _j)$ such that if $f$ is a bounded function and $|f(z_j)|\leq \eta _j$ for all $j\in \NN$ then $f\equiv 0$. It is natural to request that for sequences $(\eta _j)$ having above described property, if $z_k$ is a limit point of $(z _j)$ then $\eta _k=0$. So it is interesting that the results of Danikas and Hayman hold also for the general case of $E$ being just a non-Blaschke sequence with no more constraint. Indeed, we shall prove the following

\begin{theorem}
Assume that $E$ is a non-Blaschke sequence $(z_j)$. Then there exists a sequence of positive numbers $(\eta _j)$ with
the property that
\begin{eqnarray*}
\lim _{j\rightarrow\infty}\eta _j=0,
\end{eqnarray*}
such that if $f$ is a non-zero bounded analytic funtion on $U$ then
\begin{eqnarray*}
\limsup _{j\rightarrow \infty}\frac{|f(z_j)|}{\eta _j}=\infty .
\end{eqnarray*}
\label{GenHaymanTheorem}
\end{theorem}
Compared with results in \cite{dan} and \cite{hay}, the sequence $(\eta _j)$ in our paper is more explicitly constructed.

This paper consists of 5 parts. Some preliminary results and the proofs of Theorems \ref{ConvergenceItselfTheorem} and \ref{GenHaymanTheorem} are presented in Section 2. In Section 3, we introduce some set functions and their properties. In Section 4 we prove main result Theorem \ref{GeneralizedTheorem} applicable to sets $E$ satisfying fairly general conditions and apply it to prove Theorem \ref{LRSTheorem}. In Section 5 we gives some examples of sets $E$ that illustrate some special cases of Theorem \ref{GeneralizedTheorem}.

\section{Preliminaries}
Let us first introduce some notations and results that we are going to apply. In this section, we also give proofs of Theorems \ref{ConvergenceItselfTheorem} and \ref{GenHaymanTheorem}.

Let $Z_n$ $=$ $(z_1, z_2, \ldots, z_n)$ be a sequence of $n$ distinct points in $U$. By $B(Z_n, z)$ and $B_k(Z_n,z)$ we mean the products
\begin{equation}
B(Z_n,z) = \prod_{j = 1}^n \frac{z - z_j}{1 - \overline{z_j}z},\h B_k(Z_n,z) = \prod_{1 \leq j \neq k \leq n} \frac{z - z_j}{1 - \overline{z_j}z}.
\end{equation}

The following theorem will be used.
\begin{theorem}
If $Z_n$ = $(z_1, z_2, \ldots, z_n)$ is a sequence of $n$ distinct points in $U$ then, for all $f$ in $H^p$ and $z$ in $U$, the following inequality holds:
\begin{equation}
\left|f(z) - \sum_{k = 1}^n c_k(Z_n,z)f(z_k)\right| \leq \frac{\|f\|_p}{(1 - |z|^2)^{\frac{1}{p}}}|B(Z_n,z)|,
\label{rtieqn}
\end{equation}
where
\begin{equation}
c_{p,k}(Z_n,z) = \frac{1 - |z_k|^2}{z - z_k}\left(\frac{1 - \overline{z}z_k}{1 - |z|^2}\right)^{\frac{2 - p}{p}}\frac{B(Z_n,z)}{B_k(Z_n,z_k)}.
\label{rtckexpr}
\end{equation}
\label{RecoveryTheorem}
\end{theorem}

The reader is referred to \cite{osi} or \cite{TLNT} for proof of this Theorem.

We need some estimations of $B(Z_n, z)$ and $B_k(Z_n, z)$.

\begin{proposition}
\begin{align*}
|B(Z_n, z)| &\leq \exp\left(-\frac{1 - |z|^2}{4}\sum_{j = 1}^n (1 - |z_j|)\right),\\
|B_k(Z_n, z)| &\leq 2\exp\left(-\frac{1 - |z|^2}{4}\sum_{j = 1}^n (1 - |z_j|)\right),
\end{align*}
for $z$ in $\overline U$ and $Z_n$ in $\overline U^n$.
\label{BlaschkeProductsEstimationsProp}
\end{proposition}

\begin{proof}{}
We need to estimate $\displaystyle\left|\frac{z - \zeta}{1 - \overline{\zeta}z}\right|$ for $z$ and $\zeta$ in $U$. We compute
\begin{align*}
1 - \left|\frac{z - \zeta}{1 - \overline{\zeta}z}\right|^2
    &= \frac{(1 - |z|^2)(1 - |\zeta|^2)}{|1 - \overline{\zeta}z|^2} \geq \frac{(1 - |z|^2)(1 - |\zeta|^2)}{(1 + |\zeta|)^2}\\
    &\geq \frac{(1 - |z|^2)(1 - |\zeta|)}{2}.
\end{align*}
It follows that
\[
1 - \left|\frac{z - \zeta}{1 - \overline{\zeta}z}\right| \geq \frac{(1 - |z|^2)(1 - |\zeta|)}{2\left(1 + \left|\frac{z - \zeta}{1 - \overline{\zeta}z}\right|\right)} \geq \frac{1 - |z|^2}{4}(1 - |\zeta|),
\]
and thus,
\[
\left|\frac{z - \zeta}{1 - \overline{\zeta}z}\right| \leq \exp\left(\left|\frac{z - \zeta}{1 - \overline{\zeta}z}\right| - 1\right) \leq \exp\left(-\frac{1 - |z|^2}{4}(1 - |\zeta|)\right).
\]

From this and the fact that $\displaystyle\exp\left(\frac{1 - |z|^2}{4}(1 - |z_k|)\right)$ $\leq$ $2$, we get the assertions.
\end{proof}

Now we prove Theorems \ref{ConvergenceItselfTheorem} and \ref{GenHaymanTheorem}.

\begin{proof}{of Theorem \ref{ConvergenceItselfTheorem}}
Since $E$ satisfies (B), there is a sequence $\{z_k\}$ of distinct points in $E$ such that $\displaystyle\sum_{k=1}^\infty(1 - |z_k|)$ $=$ $\infty$.

Suppose that $f$ is a function in $H^{p}$ such that $\|f\|_1$ $\leq$ $1$ and $\displaystyle\max_{\zeta \in E}f(\zeta)$ $\leq$ $\varepsilon$. An application of Theorem \ref{RecoveryTheorem} to $Z_n$ $=$ $(z_1, z_2, \ldots, z_n)$ yields
\[
|f(z)| \leq \left|\sum _{k = 1}^n f(z_k)c_k(Z_n,z)\right|+\frac{\|f \|_1}{(1 - |z|^2)^{\frac{1}{p}}}|B(Z_n,z)|, \h\mbox{$\forall z$ $\in$ $U$}.
\]

It is clear that $c_k(Z_n, .)$ is continuous in $U$. As a consequence, we can find a constant $M_n$ depending on the sequence $\{z_k\}$ so that $|c_k(Z_n, z)|$ $\leq$ $M_n$ for all $z$ $\in$ $U$, $k$ $=$ $1, 2, \ldots, n$.

Therefore
\[
|f(z)| \leq \varepsilon n M_n + \frac{1}{(1 - |z|^2)^{\frac{1}{p}}}|B(Z_n,z)|, \h\mbox{$\forall z$ $\in$ $U$},
\]
and
\[
C_p(\varepsilon, R) \leq \varepsilon n M_n + \frac{1}{(1 - R^2)^{\frac{1}{p}}} \max_{z \in B(0,R)} |B(Z_n, z)|.
\]

Letting $\varepsilon$ $\rightarrow$ $0$ gives
\begin{equation}
\limsup_{\varepsilon \rightarrow 0} C_p(\varepsilon, R) \leq \frac{1}{(1 - R^2)^{\frac{1}{p}}} \max_{z \in B(0,R)} |B(Z_n, z)|.
\label{a}
\end{equation}

Hence to prove our assertion, it suffices to show that $\displaystyle\max_{z \in B(0, R)}B(Z_n, z)$ tends to $0$ as $n$ $\rightarrow$ $\infty$. Applying Proposition \ref{BlaschkeProductsEstimationsProp}, one get
\[
|B(Z_n,z)| \leq \exp\left(-\frac{1 - |z|^2}{4}\sum_{j = 1}^n(1 - |z_j|)\right).
\]
Thus, since $\{z_k\}$ satisfies (B), $\displaystyle\lim_{n \rightarrow \infty}\max_{z \in B(0, R)}|B(Z_n, z)|$ $=$ $0$. Hence, by passing (\ref{a}) to limit, we obtain the desired result, $\displaystyle\lim_{\varepsilon \rightarrow 0} C_p(\varepsilon, R)$ $=$ $0$.
\end{proof}

\begin{proof}{ of Theorem \ref{GenHaymanTheorem}}  From properties of $E$, we can choose a sequence of integers $n_1<n_2<...<n_k<...$ such that
\begin{eqnarray*}
\sum _{j=n_k}^{n_{k+1}-1}(1-|z_j|)\geq k.
\end{eqnarray*}

It follows that $m_k=n_{k+1}-n_{k}\geq k$. We denote $Z_{k,m_k}=\{z_{n_k},z_{n_{k}+1},...,z_{n_{k+1}-1}\}$. We define the sequence $\eta _j$ as follows
\begin{eqnarray*}
\eta _j=\frac{|B_j(Z_{k,m(k)},z_j)|}{m(k)},
\end{eqnarray*}
if $n_k\leq j<n_{k+1}$. It is easy to see that $\eta _j\rightarrow 0$ as $j\rightarrow\infty$.

Now assume that $f$ is a bounded analytic function satisfying
\begin{eqnarray*}
\limsup _{j\rightarrow \infty}\frac{|f(z_j)|}{\eta _j}<\infty ,
\end{eqnarray*}
we will show that $f\equiv 0$. Indeed, fixed $z\in U$ with $|z|\leq 1/2$. Applying Theorem \ref{RecoveryTheorem} for $Z_{k,m_{k}}$ and using Proposition \ref{BlaschkeProductsEstimationsProp} we have
\begin{eqnarray*}
|f(z)|&\leq&C(1+\sum _{j=n_{k}}^{n_{k+1}-1}\frac{|f(z_j)|}{m_k\eta _j})\max _{j=n_k,...,n_{k+1}-1}|B_j(Z_{k,m_{k}},z)|\\
&\leq&C\exp \{(-k+1)/4\},
\end{eqnarray*}
for all $k$. So $f(z)=0$ for all $|z|\leq 1/2$. Hence $f\equiv 0$.
\end{proof}
\section{Some set functions}

For results of this section forward, we need the following definitions and results in \cite{hay}.

\begin{definition}
For each set $E$ of $U$, we denote by $E_0$ the set of nontangential limit points of $E$, that is, points $\zeta$ of $\partial U$ being such that there exists a sequence $(z_n)$ in $E$ which tends nontangentially to $\zeta$, that is, such that
\begin{eqnarray*}
z_n\rightarrow \zeta,~|z_n-\zeta|=O(1-|z_n|).
\end{eqnarray*}
\end{definition}

Hayman proved the following theorem (see Theorem 1 in \cite{hay})

\noindent\textbf{Theorem A} \textit {If the set $E_0$ of nontangential limit points of a set $E$ has positive linear measure and if $f$ is a bounded analytic function satisfying
\begin{eqnarray*}
\lim _{z\in E,~|z|\rightarrow 1}f(z)=0,
\end{eqnarray*}
then $f\equiv 0$. Conversely, if $E_0$ has measure zero, then there exists $f(z)$, such that $0<|f(z)|<1$ in $U$, and satisfying
\begin{eqnarray*}
\lim _{z\in E,~|z|\rightarrow 1}f(z)=0.
\end{eqnarray*} }

Hayman's Theorem suggests us to define the following geometric condition of $E$

\begin{description}
\item[(G)] The set of nontangentially limit points $E_0$ of $E$ has measure zero.
\end{description}

From now on  we assume that $E$ satisfies condition $(G)$.

\label{setfuncs}

First, we define a new concept, the generalized Blaschke product, or the Blaschke product with respect to a function $q$. Let $q(z)$ be a function that satisfies properties of the function $f$ in Theorem A. If $\overline{E}\cap \partial U=\emptyset$, we take $q\equiv 1$. We also normalize $q$ so that $||q||_{\infty ,U}=1$. Here we use the convention that a product which has no factor is $1$. For $z$ $\in$ $\overline U$ and $Z_n$ $\in$ $\overline U^n$, we set
\begin{equation}
B_q(Z_n, z) = B(Z_n, z)q(z),\h B_{q,j}(Z_n, z) = B_j(Z_n, z)q(z).
\end{equation}

The following result is easy to prove.

\begin{proposition}
For all $R$ $\in$ $(0,1)$ and $Z_n$ $\in$ $U^n$, we have
\[
\max_{|z| \leq R} |B_q(Z_n, z)| \geq |q(0)|\prod_{1 \leq j \leq n}\max\{R, |z_j|\}.
\]
\label{GeneralizedBlaschkeProductEstimationProp}
\end{proposition}

Denote by $d(.,.)$ the Gleason distance on $U$, $d(z_1, z_2)$ $=$ $\displaystyle\left|\frac{z_1 - z_2}{1 - \overline{z_2}z_1}\right|$. Let us introduce some definitions.
\begin{definition}
For $Z_n$ $=$ $(z_1, z_2, \ldots, z_n)$ $\in$ $\overline{U}^n$, put
\begin{align}
V(Z_n)
    &= \prod_{1 \leq j \leq n}|B_q(Z_{j-1}, z_j)|,\\
\mu(z_1, z_2, \ldots, z_n)
    &= \sum_{1 \leq j \leq n} |B_j(Z_n, z_j)|^{-1},\\
M(z_1, z_2, \ldots, z_n)
    &= \sup_{z \in E}|B_q(Z_n, z)|.
\end{align}
\end{definition}
\begin{definition}
Let $E$ be a subset of $\overline{U}$ which contains infinitively many points. Put
\begin{align}
V_n(E)
    &= \sup_{Z_n \in E^n} V(Z_n),\\
\mu_n(E)
    &= \inf_{Z_n \in E^n,~V(Z_n)=V_n(E)} \mu(Z_n),\\
M_n(E)
    &= \inf_{Z_n \in E^n,~V(Z_n)=V_n(E)} M(Z_n).
\end{align}
\end{definition}

If there is no confusing, we drop argument $E$ and write $V_n,M_n,\mu _n$instead.

\begin{remark}
If $V_n$ $=$ $V(z_1, z_2, \ldots, z_n)$ then $z_j$ $\in$ $\overline E \cap U$ for all $j$ $=$ $1$, $2$, \ldots, $n$, and $|B_{q,j}(Z_n, z_j)|$ $=$ $\displaystyle \sup_{z \in E}|B_{q,j}(Z_n, z)|$.
\label{TildeV_TildeMRemark}
\end{remark}

\begin{proposition}
Let $z_1$, $z_2$, \ldots, $z_n$ and $\zeta_1$, $\zeta_2$, \ldots, $\zeta_{n+1}$ be points in $\overline E$ such that $ V(z_1, z_2, \ldots, z_n)$ $=$ $ V_n$ and $ V(\zeta_1, \zeta_2, \ldots, \zeta_{n+1})$ $=$ $ V_{n+1}$, then
\[
\mu(\zeta_1, \zeta_2, \ldots, \zeta_{n+1}) M(z_1, z_2, \ldots, z_n) \leq (n+1)\|q\|_{\infty, E}.
\]
\label{Mu_TildeMRelationProp}
\end{proposition}

\begin{proof}{}
If $z_0$ is the point in $\overline E$ such that $|B_q(Z_n, z_0)|$ $=$ $\displaystyle\prod_{1 \leq j \leq n}d(z_0, z_j) |q(z)|$ $=$ $ M(z_1, z_2, \ldots, z_n)$, then $ M(z_1, z_2, \ldots, z_n)  V(z_1, z_2, \ldots, z_n)$ $=$ $ V(z_0, z_1, z_2, \ldots, z_n)$ $\leq$ $ V_{n+1}$.

Therefore, for $k$ $=$ $1$, $2$, \ldots, $n+1$, we have
\begin{align*}
 M(z_1, z_2, \ldots, z_n)
    &\leq \frac{ V_{n+1}}{ V(z_1, z_2, \ldots, z_n)} \leq \frac{ V(\zeta_1, \zeta_2, \ldots, \zeta_{n+1})}{ V(\zeta_1, \ldots, \zeta_{k-1}, \zeta_{k+1}, \ldots, \zeta_n)}\\
    &= \prod_{1 \leq j \neq k \leq n+1}d(\zeta_k, \zeta_j)q(\zeta_k) \leq \|q\|_{\infty,E}\prod_{1 \leq j \neq k \leq n+1}d(\zeta_k, \zeta_j).
\end{align*}
It follows that
\[
\mu(\zeta_1, \zeta_2, \ldots, \zeta_{n+1}) \leq \frac{(n+1)\|q\|_{\infty,E}}{ M(z_1, z_2, \ldots, z_n)}.
\]
This proves the proposition.
\end{proof}

\begin{proposition}
$\displaystyle \lim_{n \rightarrow \infty}  V_n^{1/n}$ $=$ $\displaystyle \lim_{n \rightarrow \infty}  M_n$ $=$ $0$.
\label{ConvergenceOfTildeV_TildeMProp}
\end{proposition}

\begin{proof}{}
We need to consider only the case in which $\overline E$ intersect $\partial U$.

If $m$ $\neq$ $0$. Fix a $\delta$ $>$ $0$. By properties of $q(z)$, it follows that there exist an $r_\delta$ $<$ $1$ such that $|z|$ $<$ $r_\delta$ whenever $z$ $\in$ $\overline E$ and $q(z)$ $>$ $\delta$. For each n, we rearrange $z_1$, $z_2$, \ldots, $z_n$ so that the first $k_n$ elements have norm at least $r_\delta$ and the rest are not. We have
\begin{align*}
 V_n
    &= \prod_{1 \leq j < l \leq n}d(z_j, z_l) \prod_{1 \leq j \leq n}|q(z_j)|\\
    &\leq \|q\|_{\infty, U}^{n - k_n}\prod_{k_n + 1 \leq j < l \leq n}d(z_j, z_l) \prod_{1 \leq j \leq k_n}|q(z_j)| \leq \|q\|_{\infty, U}^n \eta^{(n - k_n)(n - k_n - 1)/2} \delta^{k_n},
\end{align*}
where $\eta$ $=$ $\frac{2r_\delta}{1 + r_\delta^2}$. It follows that $ V_n^{1/n}$ $\leq$ $\|q\|_{\infty, U}\eta^{(n - k_n)(n - k_n - 1)/2n}\delta^{k_n/n}$. From this, we see that, if $k_n/n$ $\geq$ $1/3$, then $ V_n^{1/n}$ $\leq$ $\|q\|_{\infty, U}\delta^{1/3}$, and if  $k_n/n$ $<$ $1/3$, then $ V_n^{1/n}$ $\leq$ $\|q\|_{\infty, U}\eta^{n/9}$. Hence
\[
\limsup_{n \rightarrow \infty} V_n^{1/n} \leq \limsup_{n \rightarrow \infty} \max\{\|q\|_{\infty, U}\delta^{1/3}, \|q\|_{\infty, U}\eta^{n/9}\} = \|q\|_{\infty, U}\delta^{1/3}.
\]
Since $\delta$ can be chosen arbitrarily, we deduce $\displaystyle\lim_{n \rightarrow \infty} V_n^{1/n}$ $=$ $0$.

To prove the second part of (iii), we choose $z_1$, $z_2$, \ldots, $z_n$ so that $ M_n$ $=$ $ M(z_1, z_2, \ldots, z_n)$. Noting that $|B_q(Z_n,z)|$ $\leq$ $|B_{q,j}(Z_n,z)|$ and $|B_{q,j}(Z_n,z_j)|$ $\leq$ $|B_{q,j}(Z_{j-1},z_j)|$ for all $j$ = $1$, $2$, \ldots, $n$, we compute
\begin{align*}
 M_n
    &=  M(z_1, z_2, \ldots, z_n) = \sup_{z \in E}|B_q(Z_n,z)| \leq \left(\prod_{1 \leq j \leq n} \sup_{z \in E}|B_{q,j}(Z_n, z)|\right)^{1/n}\\
    &\leq \left(\prod_{1 \leq j \leq n} |B_{q,j}(Z_n, z_j)|\right)^{1/n} \leq \left(\prod_{1 \leq j \leq n} |B_q(Z_{j-1}, z_j)|\right)^{1/n} =  V_n^{1/n}.
\end{align*}
This leads to the convergence of $M_n$ to $0$.
\end{proof}

Now we assume following condition imposed on $E$: There exists a continuous function $h:[1,\infty )\rightarrow (0,\infty )$ such that $h$ is  non-increasing, $\lim _{x\rightarrow\infty}h(x)=0$ and $M_n\leq h(n)$for all $n\in\NN$. We can define such an $h$ as follows: First, define $h(n)=\sup _{k\geq n}M_k$. Then $h(n+1)\leq h(n)$, and by Lemma \ref{ConvergenceOfTildeV_TildeMProp}, we see that $\lim _{n\rightarrow\infty} h(n)=0$. Then we extend it appropriately.

We take $\epsilon _0=\displaystyle{\frac{h(1)}{2}}$. Since $\displaystyle{\frac{h(x)}{x+1}}$ is continuous and strictly decreasing, and $\lim _{x\rightarrow\infty}h(x)=0$, we can define a function $\varphi :(0,\epsilon _0)\rightarrow (0,\infty)$ as follows:
\begin{eqnarray*}
\varphi (\epsilon )=h(x)\mbox{ iff }\epsilon =\frac{h(x)}{x+1}.
\end{eqnarray*}

We note that $\varphi $ is non-decreasing and $\lim _{\epsilon \rightarrow 0}\varphi (\epsilon)=0$.

The following result is used for proving results in Section 4

\begin{lemma}
If $R$ is a real number in $(0,1)$, then there exists a positive number $\alpha$ depending on $R$ such that for all $r$ in $[0,1]$, the inequality underneath holds,
\begin{equation}
\max\{R^\alpha, r^\alpha\} \geq \frac{R + r}{1 + Rr}.
\label{Rr}
\end{equation}
\label{ExistenceOfAlphaLemma}
\end{lemma}

\begin{proof}{}
First, we consider the case $r$ $\leq$ $R$. We have $\max\{R^\alpha, r^\alpha\}$ $=$ $R^\alpha$ and $\displaystyle\frac{R + r}{1 + Rr}$ $\leq$ $\displaystyle\frac{2R}{1 + R^2}$. Thus, if this is the case, we must choose $\alpha$ in such a way that $0$ $<$ $\alpha$ $\leq$ $\displaystyle\frac{\ln(2R) - \ln(1 + R^2)}{\ln R}$.

Finally, we consider the case $r$ $>$ $R$. The inequality (\ref{Rr}) is now equivalent to $\displaystyle\frac{r^\alpha - r}{1 - r^{\alpha + 1}}$ $\geq$ $R$. We will show that the function $f(r)$ $=$ $\displaystyle\frac{r^\alpha - r}{1 - r^{\alpha + 1}}$ ($R$ $\leq$ $r$ $\leq$ $1$) attains its absolute minimum at $R$. We have $f'(r)$ $=$ $\displaystyle\frac{r^{2\alpha} - \alpha r^{\alpha + 1} + \alpha r^{\alpha - 1} - 1}{(1 - r^{\alpha + 1})^2}$. Put $g(r)$ $=$ $r^{2\alpha} - \alpha r^{\alpha + 1} + \alpha r^{\alpha - 1} - 1$ ($R$ $\leq$ $r$ $\leq$ $1$), then $g'(r)$ $=$ $2\alpha r^{2\alpha - 1} - \alpha(1 + \alpha)r^\alpha - \alpha(1 - \alpha)r^{\alpha - 2}$. By Holder inequality, one has $(1 + \alpha)r^{1 - \alpha} + \alpha(1 - \alpha)r^{-(1 + \alpha)}$ $\geq$ $2r^{(1 - \alpha)(1 + \alpha)}r^{-(1 + \alpha)(1 - \alpha)}$ $=$ $2$ if $0$ $<$ $r$ $<$ $1$. This shows that $g(r)$ $\leq$ $0$ and thus $g(r)$ $\geq$ $g(1)$ $=$ $0$. As a consequence, $f(r)$ $\geq$ $f(R)$ $=$ $\displaystyle\frac{R^\alpha - R}{1 - R^{\alpha + 1}}$.

Therefore, the proof of the lemma is complete once we can show that for sufficiently small $\alpha$ the inequality $f(R)$ $\geq$ $R$ holds. Indeed, this is equivalent to $R^{\alpha - 1} + R^{\alpha + 1}$ $\geq$ $2$. Since $0$ $<$ $R$ $<$ $1$, it follows that $R^{-1} + R^1$ $>$ $2$. Hence, choosing $\alpha$ small enough will lead to the desired result.
\end{proof}

\section{Main result}

\begin{theorem}
Let $E$ be a subset of $U$ satifying condition $(G)$. Let $\varphi$ be the function defined in last part of Section 3. Then, there exists
a positive number $\varepsilon_0$ and a non-decreasing function $g_R:$ $(0,1]$ $\rightarrow$ $\RR^+$ satisfying
\begin{enumerate}[(i)]
\item $\displaystyle\lim_{\varepsilon \rightarrow 0}g_R(\varepsilon)$ $=$ $0$ if $E$ satisfies (B).

\item $(g_R(\varepsilon))^{[k]+1}$ $\leq$ $g_R(\varepsilon^\kappa)$ $\leq$ $(g_R(\varepsilon))^{\frac{1}{[1/\kappa]+1}}$ for positive $\kappa$.
\end{enumerate}
such that
\begin{equation}
g_R(\varepsilon) \leq C_p(\varepsilon, R) \leq Cg_R(\varphi (\varepsilon)) \h \forall \varepsilon \in (0,\varepsilon_0).
\end{equation}

In particular, if there are constants $C$ and $N$ and $\sigma >0$ such that $M_n(E)$ $\leq$ $Cn^{-\sigma}$ for $n$ $\geq$ $N$ then
\begin{equation}
g_R(\varepsilon) \leq C_p(\varepsilon, R) \leq Cg_R(\varepsilon)^{\frac{\alpha}{[1/\kappa_0]+1}} \h \forall \varepsilon \in (0,\varepsilon_0),
\end{equation}
where $\alpha$ depends on $R$ and $\sigma$.
\label{GeneralizedTheorem}
\end{theorem}

\begin{proof}{of Theorem \ref{GeneralizedTheorem} (Part 1: Existence of $g$)}
In this part, we prove the existence of $g_R$ and show that it satisfies the inequality $g_R(\varepsilon)$ $\leq$ $C_p(\varepsilon, R)$ and the two properties (i), (ii).

For $R$ in $(0,1)$, we construct the function $g_R:$ $(0,1]$ $\rightarrow$ $\RR^+$ as follows
\begin{equation}
g_R(\varepsilon) = \sup\left\{\frac{\displaystyle\prod_{1 \leq j \leq n}\max\{R^s, |z_j|^s\}}{\|q\|_{\infty, U}^s}: s,n \in \NN, Z_n \in E^n, \frac{ M^s(Z_n)}{\|q\|_{\infty, U}^s} \leq \varepsilon\right\},
\label{gRDefinitioneqn}
\end{equation}
where $ M$ is defined as in Section \ref{setfuncs}.

It can be derived from Proposition \ref{BlaschkeProductsEstimationsProp} that $\displaystyle\frac{ M^s(Z_n)}{\|q\|_{\infty, U}^s}$ $\leq$ $\displaystyle\exp\left(-s(1 - |z_0|^2)\sum_{k = 1}^n (1 - |z_k|)\right)$ for some $z_0$ in $\overline E$. Thus $g_R$ is well-defined.

\begin{enumerate}[1.]
\item $g_R(\varepsilon)$ $\leq$ $C_p(\varepsilon, R)$.

Let $z_1$, $z_2$, \ldots, $z_n$ be in $\overline E$ such that $ M^s(Z_n)$ $\leq$ $\varepsilon\|q\|_{\infty, U}^s$, then the function $h(z)$ $=$ $\displaystyle\frac{1}{\|q\|_{\infty, U}^s}B_q^s(Z_n, z)$ satisfies $|h(z)|$ $\leq$ $\varepsilon$ for $z$ $\in$ $E$ and $|h(z)|$ $\leq$ $1$ for $z$ $\in$ $U$. Hence $\displaystyle\sup_{|z| \leq R}|h(z)|$ $\leq$ $C_p(\varepsilon, R)$. On the other hand, in view of Proposition \ref{GeneralizedBlaschkeProductEstimationProp}, we have $\displaystyle\sup_{|z| \leq R}|h(z)|$ $\geq$ $\displaystyle\frac{1}{\|q\|_{\infty, U}^s}\prod_{1 \leq j \leq n}\max\{R^s, |z_j|^s\}$. This implies that $g_R(\varepsilon)$ $\leq$ $C_p(\varepsilon,R)$.

\item $\displaystyle\lim_{\varepsilon \rightarrow 0} g_R(\varepsilon)$ $=$ $0$ if $E$ satisfies (B).

This assertion follows from 1. and Theorem \ref{ConvergenceItselfTheorem}.

\item $(g_R(\varepsilon))^{[k]+1}$ $\leq$ $g_R(\varepsilon^\kappa)$ $\leq$ $(g_R(\varepsilon))^{\frac{1}{[1/\kappa]+1}}$ for positive $\kappa$.

It is clear that, for arbitrary $\varepsilon$ in $(0, 1)$, $ M^s(Z_n)$ $\leq$ $\|q\|_{\infty,U}^s\varepsilon^\kappa$ implies immediately that $ M^{s([1/\kappa]+1)}(Z_n)$ $\leq$ $\|q\|_{\infty, U}^{s([1/\kappa]+1)}\varepsilon$. Therefore, $g_R(\varepsilon^\kappa)$ $\leq$ $(g_R(\varepsilon))^{\frac{1}{[1/\kappa]+1}}$.

Similarly, we get $g_R(\varepsilon)$ $\leq$ $(g_R(\varepsilon^\kappa))^{\frac{1}{[\kappa]+1}}$. Hence $g_R(\varepsilon^\kappa)$ $\geq$ $(g_R(\varepsilon))^{[\kappa]+1}$.
\end{enumerate}
This ends part 1.
\end{proof}

\begin{remark}
For two different chooses $q(z)=q_1(z)$ and $q(z)=q_2(z)$, it is not expected that $g_1(\epsilon )=g_2(\epsilon )$. But $g_1(\epsilon )$ and $g_2(\epsilon )$are still related to each other, since according to Theorem \ref{GeneralizedTheorem} we have
\begin{eqnarray*}
g_1(\epsilon )\leq C_p(\epsilon )\leq Cg_2(\varphi _2(\epsilon)).
\end{eqnarray*}

Here we droped indices $R$ in functions $g_1,~g_2$ and $C_p$.

In particular, if there are constants $C$ and $N$ and $\sigma >0$ such that $M_n(E)$ $\leq$ $Cn^{-\sigma}$ for $n$ $\geq$ $N$ then we can bound $g_1$ by a positive power of $g_2$.
\end{remark}

Theorem \ref{RecoveryTheorem} provides a way to investigate $C_p(\varepsilon, R)$, that is to estimate $D_p(\varepsilon,R)$ which will be defined underneath.
\begin{definition}
For $Z_n$ $=$ $(z_1$, $z_2$, \ldots, $z_n)$ $\in$ $E^n$, put
\begin{equation}
D_p(\varepsilon, R, Z_n) = \sup_{|z| \leq R}\left(\varepsilon\sum_{k = 1}^n|c_k(Z_n, z)| + \frac{1}{(1 - |z|^2)^{\frac{1}{p}}}|B(Z_n,z)|\right).
\end{equation}
The infimum of $D_p(\varepsilon, R, Z_n)$ while $Z_n$ varies in $E^n$ is defined to be $D_{p,n}(\varepsilon, R)$,
\begin{equation}
D_{p,n}(\varepsilon, R) = \inf_{Z_n \in E^n}D(\varepsilon, R, Z_n).
\end{equation}
Finally, set
\begin{equation}
D_p(\varepsilon, R) = \inf_{n \in \NN}D_{p,n}(\varepsilon, R).
\end{equation}
\end{definition}

Recall that
\[
c_{p,k}(Z_n,z) = \frac{1 - |z_k|^2}{z - z_k}\left(\frac{1 - \overline{z}z_k}{1 - |z|^2}\right)^{\frac{2 - p}{p}}\frac{B(Z_n,z)}{B_k(Z_n,z_k)}.
\]

We have
\[
D_p(\varepsilon, R, Z_n) \leq C E(z_1, z_2, \ldots, z_n)(\varepsilon\mu(z_1, z_2, \ldots, z_n) + 1),
\]
where $\displaystyle E(z_1, z_2, \ldots, z_n)$ $=$ $\displaystyle\prod_{k = 1}^n\frac{R + |z_k|}{1 + |z_k|R}$.

\begin{proof}{of Theorem \ref{GeneralizedTheorem} (Part 2: finished)}
It follows from Lemma \ref{ExistenceOfAlphaLemma} that
\begin{equation}
D_{p,n}(\varepsilon, R) \leq C \prod_{1 \leq j \leq n}\max\{R^\alpha, |\zeta_j|^\alpha\}(\varepsilon\mu(\zeta_1, \zeta_2, \ldots, \zeta_n)+1),
\end{equation}
where $\zeta_1$, $\zeta_2$, \ldots, $\zeta_n$ are defined by $ V(\zeta_1, \zeta_2, \ldots, \zeta_n)$ $=$ $ V_n(E)$ and $M_n(\zeta_1, \zeta_2, \ldots, \zeta_n)=M_n(E)$. Note that these $\zeta_j$ is contained in $U$. For simplicity, we denote by $\mu_n^*$ the quantity $\mu(\zeta_1, \zeta_2, \ldots, \zeta_n)$.

It follows from Proposition \ref{ConvergenceOfTildeV_TildeMProp} that $\displaystyle\lim_{n \rightarrow \infty}  M_n(E)$ $=$ $0$. Thus, we can choose the smallest $n_0$ such that $ M_{n_0}(E)$ $\leq$ $\varphi (\varepsilon )$ $<$ $ M_{n_0-1}(E)$ for all $\varepsilon$ less than some fixed constant $\varepsilon_0$. Then, by Proposition \ref{Mu_TildeMRelationProp}
\[
 M_{n_0}(E) \leq \varphi (\varepsilon ) <  M_{n_0-1}(E) \leq \frac{n_0\|q\|_{\infty,E}}{\mu_{n_0}^*}.
\]

On the other hand, we have $\varphi (\epsilon ) < M_{n_0-1}(E)$ $\leq$ $h(n_0-1)$ for $n_0$ $\geq$ $N$. This gives $n_0$ $\leq$ $x+1$, where
\begin{eqnarray*}
\epsilon =\frac{h(x)}{x+1}.
\end{eqnarray*}
 Hence,
\begin{eqnarray*}
\epsilon \mu _{n_0}^{*}=\frac{\epsilon}{\varphi (\epsilon )}\varphi (\epsilon )\mu _{n_0}^{*}\leq\frac{\epsilon}{\varphi (\epsilon )}n_0
\leq\frac{\epsilon}{\varphi (\epsilon )}(x+1)=\frac{\frac{h(x)}{x+1}}{h(x)}(x+1)=1.
\end{eqnarray*}

Now, noting that $\|q\|_{\infty, U}$ $\geq$ $1$, it is deduced that $ M(\zeta_1, \zeta_2, \ldots, \zeta_{n_0})$ $=$ $ M_{n_0}(\overline E)$ $\leq$ $\varphi (\varepsilon )$ $\leq$ $\varphi (\varepsilon )\|q\|_{\infty, U}$. This yields
\[
\prod_{1 \leq j \leq n}\max\{R, |\zeta_j|\} \leq \|q\|_{\infty, U}g_R(\varphi (\varepsilon )).
\]
Thus, the assertion of the theorem follows.

To prove for the case $M_n\leq Cn^{-\sigma}$ we see that $h(x)=Cx^{-\sigma}$. So we have
\begin{eqnarray*}
\varphi (\epsilon)=h(x)=Cx^{-\sigma}\leq C_2\epsilon ^{\sigma /(1+\sigma )},
\end{eqnarray*}
since $\epsilon =Cx^{-\sigma}(1+x)^{-1}$. Now apply property (ii) of $g$ we get the result.

\end{proof}

Now we are ready to prove Theorem \ref{LRSTheorem}.

\begin{proof}{ of Theorem \ref{LRSTheorem}} Since $\overline E\subset U$ we can take $q\equiv 1$. We can assume that $E$ is closed. From the assumption, we get that $M_n(E)\leq \displaystyle{(\frac{2r}{1+r})^n}$. So the function $\varphi (\epsilon )$ in Theorem \ref{GeneralizedTheorem} satisfies
\begin{equation}
\lim _{\epsilon\rightarrow 0}\frac{\log \epsilon}{\log \varphi (\epsilon )}=1.\label{LRS.2}
\end{equation}

If we put
\begin{eqnarray*}
\widetilde {g}_p(\epsilon ,R)=\sup \{\max _{|z|\leq R}|B(Z_n,z)|:~n\in \NN,~Z_n\in E^n,~\max _{z\in E }|B(Z_n,z)|\leq \epsilon\},
\end{eqnarray*}
then it follows from the proof of Theorem \ref{GeneralizedTheorem} that there exist $C ,~\epsilon _0>0$ such that for all $0<\epsilon <\epsilon _0$ we have
\begin{equation}
\widetilde {g}_p(\epsilon ,R)\leq C_p(\epsilon ,R)\leq C\widetilde {g}_p(\epsilon ,R)^{1/(1+[\log (\varphi (\epsilon ))/\log \epsilon])}.\label{LRS.1}
\end{equation}
Indeed, to prove the left-hand side of this inequality we can use the definition of $C_p(\epsilon ,R)$, while for the right-hand side we use Theorem \ref{GeneralizedTheorem} and the property that
\begin{eqnarray*}
g_p(\epsilon ,R)\leq \widetilde{g_p}(\epsilon ,R).
\end{eqnarray*}

From (\ref{LRS.2}) we can write (\ref{LRS.1}) such as
\begin{eqnarray*}
\widetilde {g}_p(\epsilon ,R)\leq C_p(\epsilon ,R)\leq C\widetilde {g}_p(\epsilon ,R)^{1/2}.
\end{eqnarray*}

Since $E$ is compact, there exist $n\in \NN$ and $Z_n\in E^n$ such that
\begin{eqnarray*}
\widetilde {g}_p(\epsilon ,R)=\max _{|z|\leq R}|B(Z_n,z)|.
\end{eqnarray*}
So we get the conclusion of Theorem \ref{LRSTheorem}.
\end{proof}
\section{Some examples}

In this Section we gives some examples of sets satisfying the special case $M_n\leq Cn^{-\sigma}$ of Theorem \ref{GeneralizedTheorem}. The criteria given in this Section are fairly easy to verify.

\begin{proposition}
Assume that $E$ is contained in some Stolz angles. Then there exist $\sigma$, $C$ and $N$ $>$ $0$ such that $ M_n(E)$ $\leq$ $Cn^{-\sigma}$ for $n$ $\geq$ $N$.
\label{CorrectnessOfGeneralizationProp}
\end{proposition}

\begin{proof}{}

Let $\overline{E}\cap \partial U=\{a_1,a_2,...,a_n\}$. We take in this case $q(z)=(z-a_1)(z-a_2)...(z-a_n)$.

We seperate the proof into three steps.
\begin{enumerate}[1.]
\item Suppose that $\overline E$ lies inside U. In this case $M_n$ $\leq$ $n^{-\sigma}$ for sufficiently large $n$.

\item Suppose that $\overline E$ $\cap$ $\partial U$ has only one point. By mean of some rotation, we may assume that it this point is $1$.

We have $q(z)$ $=$ $z - 1$. We see that if $|q(z)|$ $>$ $\delta$ $>$ $0$ for some $z$ in $E$, then $|z|$ $<$ $r_\delta$ = $1 - c\delta$ where $c$ is a constant depending on the Stolz angle with vertex at $1$. Refering to the proof of Proposition \ref{ConvergenceOfTildeV_TildeMProp}, we get
\begin{equation}
 V_n^{1/n} \leq C\max\{\delta^{1/3}, \eta^{n/9}\}.
\label{eqmax}
\end{equation}
Choosing $\delta$ $=$ $n^{-3\sigma}$ ($\sigma$ $\in$ $(0,1/6)$), we have
\[
\eta = \frac{2r_\delta}{1 + r_\delta^2} = \frac{2(1 - c n^{-3\sigma})}{1 + (1 - c n^{-3\sigma})^2} = \frac{2n^{6\sigma} - 2cn^{3\sigma}}{2n^{6\sigma} - 2cn^{3\sigma} + c^2}.
\]
Hence,
\begin{align*}
\eta^{n/9}
    &= \left(\frac{2n^{6\sigma} - 2cn^{3\sigma}}{2n^{6\sigma} - 2cn^{3\sigma} + c^2}\right)^{n/9} = \left(1 - \frac{c^2}{2n^{6\sigma} - 2cn^{3\sigma} + c^2}\right)^{n/9}\\
    &\leq \left(1 - \frac{c^2}{2n^{6\sigma}}\right)^{\frac{2n^{6\sigma}}{c^2}\frac{c^2n^{1-6\sigma}}{18}} \leq \exp\left(-\frac{c^2n^{1-6\sigma}}{18}\right) \leq n^{-\sigma}
\end{align*}
for sufficiently large $n$. Combining with (\ref{eqmax}), the assertion follows.

\item Now, consider the general case. It suffices to show that if $E_1$ and $E_2$ are two sets satisfy $ V_n^{1/n}(E_i)$ $\leq$ $C n^{-\sigma_i}$ for $n$ $\geq$ $N$ ($i$ $=$ $1$, $2$) and $E$ $=$ $E_1$ $\cup$ $E_2$, then $ V_n^{1/n}(E)$ $\leq$ $C n^{-\sigma}$, for $n$ $\geq$ $2N$ and $\sigma$ $=$ $\min\{\sigma_1, \sigma_2\}/2$. We take $q(z)=q_1(z)q_2(z)$where $q_1,~q_2$ are coressponding $q's$ functions of $E_1,~E_2$.  Fix an $n$ $\geq$ $2N$ and suppose that $ V_n(E)$ $=$ $ V(z_1, z_2, \ldots, z_l, \zeta_1, \zeta_2, \ldots, \zeta_k)$ for $z_j$ $\in$ $E_1$, $\zeta_j$ $\in$ $E_2$, and $n$ $=$ $l + k$. It follows from definitions that
\[
 V_n^{1/n}(E) \leq  V_l^{1/n}(E_1)  V_k^{1/n}(E_2).
\]
We may assume that $l$ $\geq$ $k$. It implies that $l$ $\geq$ $n/2$ $\geq$ $N$. If $k$ $\leq$ $N$, we have
\[
 V_n^{1/n}(E) \leq C V_l^{1/n}(E_1) \leq C l^{-\sigma_1 l/n} \leq C (n/2)^{- \sigma_1/2} \leq C n^{-\sigma}.
\]
If $k$ $\geq$ $N$, we have
\begin{align*}
 V_n^{1/n}(E)
    &\leq  V_l^{1/n}(E_1)  V_k^{1/n}(E_2) \leq C l^{-\sigma_1 l/n}k^{-\sigma_2 k/n} \leq C \left(l^{-l/n}k^{-k/n}\right)^{-\sigma}\\
    &= C n^{-\sigma}\left((l/n)^{-l/n}(k/n)^{-k/n}\right)^{-\sigma} \leq C n^{-\sigma}.
\end{align*}
Here we have used the inequality $x^x(1-x)^{1-x}$ $\geq$ $1/2$ for all $x$ $\in$ $(0,1)$.
\end{enumerate}
The proof is complete.
\end{proof}

To find a more sophisticated example, we use the construction used in \cite{hay}. For convenience, we recall some definitions that Hayman used in constructing the function $f$ in Theorem A.

\begin{definition} Let $E$ satisfy $(G)$. We write
\begin{eqnarray*}
E'=\{z=re^{i\theta}:~|\theta -\phi |< 1-r\mbox{ and }re^{i\phi} \in E \}.
\end{eqnarray*}

Next, for $0\leq\theta\leq 2\pi$, we define
\begin{eqnarray*}
\rho (\theta )=\sup \{\rho :~0\leq \rho <1,~\rho e^{i\theta}\in E'\}.
\end{eqnarray*}

Let $E_{\infty}$ be the set of $\theta$ such that $\rho (\theta )=1$. If $\theta \in E_\infty$ then $e^{i\theta}\in E_0$. So $m(E_\infty )=0$, where $m(.)$ is the Lebesgue's measure of the unit circle.

For each $1>r>0$ let $E_r$ be the set of all $\theta$ such that $0\leq \theta \leq 2\pi$ and $\rho (\theta )> r$. Then $E_r$ are open and contract with increasing $r$, and
\begin{eqnarray*}
\bigcap _{r}E_r=E_{\infty}.
\end{eqnarray*}
Thus
\begin{eqnarray*}
\lim _{r\rightarrow 1}m(E_r)=0.
\end{eqnarray*}
\end{definition}

Considering carefully the construction in the proof of Theorem 1 in \cite{hay} and Step 2 of the proof of Proposition \ref{CorrectnessOfGeneralizationProp} we can show that if the quantities $m(E_r)$ tend to $0$ sufficiently fast, then $M_n\leq Cn^{-\sigma}$. In particular, this claim is true if the following condition is satisfied
\begin{eqnarray*}
m(E_{\delta})\leq \frac{1}{-2\log \epsilon} \mbox{ if } \delta =1-K\sqrt{-\epsilon\log \epsilon},
\end{eqnarray*}
where $K$ is a positive constant. In fact, if this condition holds, the function $f$ is constructed in Theorem 1 in \cite{hay} will satisfy: if $|f(z)|>\epsilon$ then $|z|\leq 1-K\sqrt{-\epsilon\log \epsilon}$.

\end{document}